\documentclass[12]{article}
\usepackage{fleqn}
\usepackage{graphicx}

\usepackage{amsmath}
\usepackage{amssymb}
\begin{document}
\Large
\begin{center}
{\bf Veldkamp-Space Aspects of a Sequence of Nested\\ Binary Segre Varieties}
\end{center}
\vspace*{-.3cm}
\large
\begin{center}
Metod Saniga,$^{1,2}$ Hans Havlicek,$^{1}$  Fr\'ed\'eric Holweck,$^{3}$ Michel Planat$^{4}$ \\ and Petr Pracna$^{5}$ 
\end{center}
\vspace*{-.6cm}
\normalsize
\begin{center}
$^{1}$Institute for Discrete Mathematics and Geometry,
Vienna University of Technology\\ Wiedner Hauptstra\ss e 8--10,
A-1040 Vienna, Austria\\
(havlicek@geometrie.tuwien.ac.at)

\vspace*{.1cm}

$^{2}$Astronomical Institute, Slovak Academy of Sciences\\
SK-05960 Tatransk\' a Lomnica, Slovak Republic\\
(msaniga@astro.sk) 

\vspace*{.1cm}

$^{3}$Laboratoire IRTES/M3M, Universit\'e de Technologie de Belfort-Montb\'eliard\\ 
F-90010 Belfort, France\\ (frederic.holweck@utbm.fr) 

\vspace*{.1cm}

$^{4}$Institut FEMTO-ST, CNRS, 32 Avenue de l'Observatoire,
F-25044 Besan\c con, France\\
(michel.planat@femto-st.fr)

\vspace*{.1cm}

$^{5}$J. Heyrovsk\' y Institute of Physical Chemistry, v.v.i.,  Academy of Sciences of Czech Republic Dolej\v skova 3, CZ-18223 Prague, Czech Republic\\
(pracna@jh-inst.cas.cz)

\end{center}

\vspace*{-.4cm} \noindent \hrulefill

\vspace*{-.0cm} \noindent {\bf Abstract}

\noindent Let $S_{(N)} \equiv PG(1,\,2) \times PG(1,\,2) \times \cdots \times PG(1,\,2)$ be a Segre variety that is $N$-fold direct product of projective lines of size three. Given two geometric hyperplanes $H'$ and $H''$ of $S_{(N)}$,  let us call the triple $\{H', H'', \overline{H' \Delta H''}\}$ the Veldkamp line of $S_{(N)}$. We shall demonstrate, for the sequence $2 \leq N \leq 4$, that the properties of  geometric hyperplanes of $S_{(N)}$ are fully encoded in the properties of Veldkamp {\it lines} of $S_{(N-1)}$. Using this property, a complete classification of all types of geometric hyperplanes of $S_{(4)}$ is provided.  Employing the fact that, for $2 \leq N \leq 4$, the (ordinary part of) Veldkamp space of $S_{(N)}$ is $PG(2^N-1,2)$, we shall further describe which types of geometric hyperplanes of $S_{(N)}$ lie on a certain hyperbolic quadric $\mathcal{Q}_0^+(2^N-1,2) \subset PG(2^N-1,2)$ that contains the $S_{(N)}$  and is invariant under its stabilizer group; in the $N=4$ case we shall also single out those of them that correspond, via the Lagrangian Grassmannian of type $LG(4,8)$, to the set of 2295 maximal subspaces of the symplectic polar space $\mathcal{W}(7,2)$.

\vspace*{.1cm}

\noindent
{\bf Keywords:} Binary Segre Varietes -- Veldkamp Spaces -- Hyperbolic Quadrics 

\vspace*{-.2cm} \noindent \hrulefill

\vspace*{-.4cm}

\section{Introduction}
\vspace*{-.3cm}
Algebraic varieties that are Cartesian products of two or more projective spaces of the same or different dimensions were introduced as early as 1891 by Corrado Segre \cite{cs} who, however, considered only the `classical' case, viz. the field of real or complex numbers.  Some seventy years later, Benjamino Segre \cite{bs} showed that the definition of Segre varieties, as well as a majority of their properties, carry over to other fields, in particular over finite (Galois) fields (see also \cite{mo}). The past decade witnessed an increased interest in the topic, and this not only from a pure mathematical point of view (e.\,g., \cite{bcr,lan,lw,lash,cool,shaw,sgh}), but also due to important quantum physical applications. 
The latter case is remarkable in that alongside classical Segre varieties, and associated Segre maps, which provide
a setup for geometrical construction of concurrence, one of principal measures of quantum entanglement (e.\,g., \cite{bh,miy,hb,hey,HLT1}),
we also encounter Segre varieties defined over the smallest Galois field, which are intimately linked with the issue of quantum contextuality and the so-called black-hole--qubit correspondence through the properties of certain generalized multi-qubit Pauli groups (e.\,g., \cite{ps,lsvp,hsl}).
Interestingly enough, these `quantum contextual' finite Segre varieties are found to lie within the family defined as $S_{(N)} \equiv PG(1,2) \times PG(1,2) \times \cdots \times PG(1,2)$ ($N$ times), where $PG(1,2)$ represents the smallest projective line, being always contained/embedded in a distinguished finite geometry/point-line incidence structure. 

The first physically relevant type of
finite Segre variety is $S_{(2)}$. It lives in the symplectic polar space $\mathcal{W}(3,2)$, which underlies the commutation relations between the elements of two-qubit generalized Pauli group \cite{ps,hos}. The $S_{(2)}$ itself, being isomorphic to the smallest slim generalized quadrangle $GQ(2,1)$, is the geometry behind what is known as a Mermin magic square \cite{mer}. It represents a set of nine observables placed at the vertices of a $3 \times 3$ grid and forming six maximum sets of pairwise commuting elements that lie along three horizontal and three vertical lines, each observable thus pertaining to two such sets. The observables are selected in such a way that the product of their triples in five of the six sets is $+I$, whilst in the remaining set it is $-I$, $I$ being the identity matrix. 
Another prominent type of the Segre variety is $S_{(3)}$, which enters the quantum informational game in a more refined disguise --- namely through its dual $S_{(3)}^{\star}$. The latter lives in the generalized quadrangle of type $GQ(2,4)$, one of the key finite geometries in the context of the  black-hole--qubit correspondence \cite{lsvp}. $GQ(2,4)$ encodes completely the entropy formula of black holes, or black strings, in certain $D=5$ supergravity theories. This entropy formula features 27 charges and 45 terms of three charges each that correspond, respectively, to 27 points and 45 lines of $GQ(2,4)$.
Take any three pairwise disjoint $GQ(2,1)$s in  $GQ(2,4)$ and remove their lines; what is left is a copy of $S_{(3)}^{\star}$. Hence, our $S_{(3)}^{\star}$ represents an interesting sub-geometry of the entropy formula that still incorporates all 27 charges, but picks up only a particular subset of 27 terms out of the totality of 45 ones.

These findings are indicative of the fact that also the next type of Segre variety in the hierarchy, $S_{(4)}$, is --- up to duality --- likely to be of relevance for quantum physics and thus deserves a closer look at. To perform such an inspection, we shall employ the notion of Veldkamp space of a point-line incidence geometry \cite{buec}.  This important concept has already been utilized to analyze the structure of $\mathcal{W}(3,2)$ \cite{spph} and higher-rank symplectic polar spaces \cite{vrle} where it helped to reveal some finer traits of the structure of multiple-qubit Pauli groups, that of $GQ(2,4)$ \cite{sglpv} where it led to a deeper understanding of distinguished truncations of the above-described black-hole entropy formula, as well as that of the Segre variety $S_{(3)}$ \cite{grsa}. Our reasoning will heavily rest on the results of the last-mentioned reference, where properties of both geometric hyperplanes and Veldkamp lines of $S_{(3)}$ were treated in much detail. 

The paper is organized as follows. Section 2 highlights basic concepts, symbols and notation to be employed. In Section 3, Subsecs. 3.1 to 3.3, we shall describe, for the sequence $2 \leq N \leq 4$, a diagrammatical recipe of how to ascertain the properties of  Veldkamp 
{\it points} of $S_{(N)}$ if we know the properties of Veldkamp {\it lines} of $S_{(N-1)}$; using this recipe, a complete classification of geometric hyperplanes of $S_{(4)}$ will be arrived at. In addition, for each $S_{(N)}$ in the above-given sequence it shall be shown which types of its geometric hyperplanes lie on the unique hyperbolic quadric $\mathcal{Q}_0^+(2^N-1,2) \subset PG(2^N-1,2)$ that contains the $S_{(N)}$  and is invariant under its stabilizer group.
Finally, Section 4 will be devoted to a brief summary of the main findings as well as to pointing out certain resemblance 
between our `generalized' concept of Veldkmap space and the notion of projective space defined over a ring.

\vspace*{-.3cm}
\section{Basic concepts, symbols and notation}
\vspace*{-.3cm}
In this section we shall make a brief inventory of basic concepts, symbols and notation employed in the sequel.

Our starting point is a {\it point-line incidence structure} $\mathcal{C} = (\mathcal{P},\mathcal{L},I)$ where $\mathcal{P}$ and $\mathcal{L}$ are, respectively, sets of points and lines and where incidence $I \subseteq \mathcal{P} \times \mathcal{L}$ is a binary relation indicating which point-line pairs are incident (see, e.\,g., \cite{shult}).  The dual of a point-line incidence structure is the structure with points
and lines exchanged, and with the reversed incidence relation. In what follows we shall encounter only specific point-line incidence structures where every line has the same number of points, every point is incident with the same number of lines and any two distinct points are joined by at most one line. A {\it geometric hyperplane} of $\mathcal{C} = (\mathcal{P},\mathcal{L},I)$ is a proper subset of $\mathcal{P}$ such that a line from $\mathcal{C}$ either lies fully in the subset, or shares with it only one point.
Given a hyperplane $H$ of $\mathcal{C}$, one defines the
{\it order} of any of its points as the number of lines through
the point that are fully contained in $H$; a point of $H$  is called {\it deep} if all the lines
passing through it are fully contained in $H$. 
 If $\mathcal{C}$ possesses geometric hyperplanes, then one can define the {\it Veldkamp space} of $\mathcal{C}$ as follows \cite{buec}: (i) a point of the Veldkamp space is a geometric hyperplane of  $\mathcal{C}$
and (ii) a line of the Veldkamp space is the collection $H'H''$ of all geometric hyperplanes $H$ of $\mathcal{C}$  such that $H' \cap H'' = H' \cap H = H'' \cap H$ or $H = H', H''$, where $H'$ and $H''$ are distinct geometric hyperplanes. There exists a wide family of point-line incidence structures, including also those to be discussed below, where each line has three points and where a line of the Veldkamp space can equivalently be defined as $\{H', H'', \overline{H' \Delta H''}\}$; here, the symbol $\Delta$ stands for the symmetric difference of the two geometric hyperplanes and an overbar denotes the complement of the object indicated.  Occasionally, $\overline{H' \Delta H''}$ will be called the {\it Veldkamp sum} of hyperplanes $H'$ and $H''$.
To meet our needs, we shall use a slightly {\it generalized} notion of the Veldkamp space of $\mathcal{C}$, $\mathcal{V}(\mathcal{C})$, which also includes $\mathcal{C}$ as the extraordinary geometric hyperplane and which features extraordinary Veldkamp lines, the latter being of type $\{H, H, \mathcal{C}\}$.  
Moreover, we shall reserve the symbol ord-$\mathcal{V}(\mathcal{C})$ for the `ordinary' part of $\mathcal{V}(\mathcal{C})$, i.\,e. that devoid of the extraordinary hyperplane and all the Veldkamp lines passing through it.  

Next, let $V(d+1,q)$, $d \geq 1$, denote a rank-$(d+1)$ vector space over the Galois field $GF(q)$, $q$ being a power of a prime. Associated with this vector space is a $d$-dimensional {\it projective space} over $GF(q)$, $PG(d,q)$, whose points, lines, planes,$\ldots$, hyperplanes are rank-one, rank-two, rank-three,$\ldots$, rank-$d$ subspaces of $V(d+1,q)$.
A {\it quadric} in $PG(d, q)$, $d \geq 1$, is the set of points whose coordinates satisfy an equation of the form $\sum_{i,j=1}^{d+1} a_{ij} x_i x_j = 0$, where at least one $a_{ij} \neq 0$. 
Up to transformations of coordinates, there is one or there are two distinct kinds of non-singular quadrics in $PG(d, q)$ according as $d$ is even or odd, namely \cite{ht}:
$\mathcal{Q}(2N,q)$, the {\it parabolic} quadric formed by all points of $PG(2N, q)$ satisfying the standard equation $x_1x_2+\cdots+x_{2N-1}x_{2N} + x_{2N+1}^{2} = 0$;
$\mathcal{Q}^{-}(2N - 1,q)$, the {\it elliptic} quadric formed by all points of $PG(2N - 1, q)$ satisfying the standard equation $f(x_1,x_2)+x_3x_4+\cdots+x_{2N-1}x_{2N} = 0$, where $f$ is irreducible over GF$(q)$; and
$\mathcal{Q}^{+}(2N - 1,q)$, the {\it hyperbolic} quadric formed by all points of $PG(2N - 1, q)$ satisfying the standard equation $x_1x_2+x_3x_4+\cdots+x_{2N-1}x_{2N} = 0$;
where $N \geq 1$. 
In this paper we shall only be concerned with hyperbolic quadrics, in particular with those of $q=2$; each such quadric is found to accommodate as many as $(2^{N-1}+1)(2^{N}-1)$ points.
Given a $PG(2N-1,q)$ that is endowed with a symplectic form, the {\it symplectic polar space} $\mathcal{W}(2N-1,q)$ in $PG(2N-1,q)$ is the space of all totally isotropic subspaces with respect to the symplectic form \cite{cam}, with its maximal totally isotropic subspaces, also
called {\it generators}, having dimension $N - 1$.  For $q=2$ this
polar space contains $|PG(2N-1, 2)| = 2^{2N} - 1 = 4^{N} - 1$
points and $(2+1)(2^2+1)\cdots(2^N+1)$ generators.

Last but not least, there comes a finite Segre variety \cite{ht} $S_{d_1, d_2,\ldots,d_N}(q) \equiv PG(d_1,q) \times PG(d_2,q) \times \ldots \times PG(d_N,q),$ whose simplest possible form, namely $S_{1, 1,\ldots,1}(2) \equiv S_{(N)}$,  will be of central interest to us. An easy count yields that $S_{(N)} = PG(1,2) \times S_{(N-1)}$, when viewed as an abstract point-line incidence structure, has $3^N$ points and $N 3^{N-1}$ lines, with three points per a line and $N$ lines through a point.

\begin{table}[pth!]
\begin{center}
\caption{The two types of geometric hyperplanes of $S_{(2)}$. } 
\vspace*{0.2cm}
{\begin{tabular}{|c|c|c|c|c|c|c|c|c|c|} \hline \hline
\multicolumn{1}{|c|}{} & \multicolumn{1}{|c|}{} & \multicolumn{1}{|c|}{}  &  \multicolumn{3}{|c|}{} & \multicolumn{2}{|c|}{}  &\multicolumn{1}{|c|}{} &\multicolumn{1}{|c|}{}   \\
\multicolumn{1}{|c|}{} & \multicolumn{1}{|c|}{} &
\multicolumn{1}{|c|}{}  &  \multicolumn{3}{|c|}{Pts of
Order} & \multicolumn{2}{|c|}{$S_{(1)}$'s Type} & \multicolumn{1}{|c|}{} & \multicolumn{1}{|c|}{} \\
 \cline{4-8}
Tp & Ps  & Ls   & ~0~   & ~1~   & 2  & ~D~  & $H_1$ & VL & Crd  \\
\hline
1  & 5   &  2   & 0   &  4  & 1  & 2  & 4     & ext  & 9    \\
\hline
2  & 3   &  0   & 3   &  0  & 0  & 0  & 6     & ord  & 6   \\
\hline \hline
\end{tabular}}
\end{center}
\end{table}

\vspace*{-.8cm}
\section{From the Veldkamp lines of \textbf{\textit{S}}$_{(N-1)}$ to the Veldkamp points of \textbf{\textit{S}}$_{(N)}$ ($2 \leq N \leq 4$)}
\subsection{\textbf{\textit{S}}$_{(2)}$ and its Veldkamp space}
The Segre variety $S_{(2)}$ is a point-line incidence structure having 9 points and 6 lines, with two lines per a point.
Using its diagrammatical representation as a $3 \times 3$ grid, one can readily check that it features two distinct kinds of geometric hyperplanes, which are depicted in Figure 1 and whose properties are summarized in Table 1. In the latter, the first column gives the type (`Tp') of a hyperplane, which is followed by the number of points (`Ps') and lines (`Ls') it contains, and the number of points of given order. The next two columns tell us about how many of 6 lines ($S_{(1)}$'s) are fully located (`D') in the hyperplane and/or share with it a single point (which is the only ordinary hyperplane of $S_{(1)}$). The VL-column lists the types of Veldkamp lines of $S_{(1)}$  we get by projecting a hyperplane of the given type into a line of the $S_{(2)}$. Finally, for each hyperplane type we give its cardinality (`Crd'). A brief inspection of Figure 1 shows that, given a point of $S_{(2)}$, a type one hyperplane consists of the point and all the points that are not at maximum 
distance\footnote{Here the term `distance' means the graph-theoretical distance between the corresponding vertices of the collinearity graph of the point-line incidence structure in question.} from it; in the language of near polygons (e.\,g., \cite{bru}), this is usually termed a {\it singular} hyperplane.
On the other hand, a type two hyperplane represents a maximum set of mutually non-collinear points, often referred to as an {\it ovoid}.

\begin{figure}[t]
\centerline{\includegraphics[width=6.7truecm,clip=]{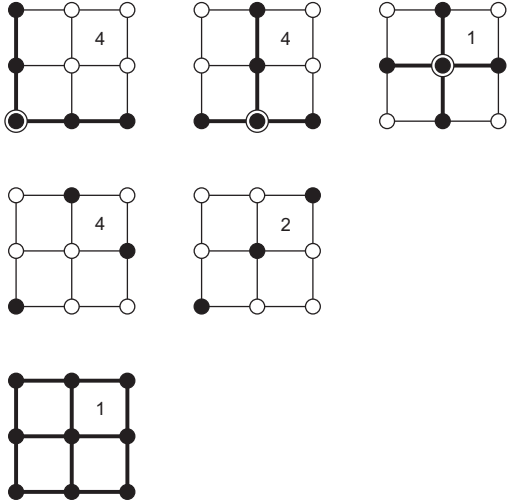}}
\caption{A diagrammatic representation of the types of geometric hyperplanes of $S_{(2)}$. The number attached to a subfigure indicates how many distinct copies of a given hyperplane one gets by rotating the subfigure around its center. The top row illustrates all 9 copies of type one hyperplane, the middle row all 6 copies of type two hyperplane and, for the sake of completeness, the bottom row shows the extraordinary hyperplane. The encircled bullets in the top row denote deep points. (Obviously, all the nine points of the extraordinary hyperplane are deep (not indicated).)}
\end{figure}
\begin{figure}[pth!]
\centerline{
\includegraphics[width=4.2truecm,clip=]{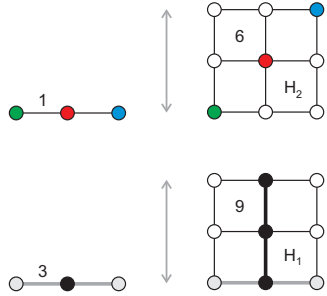}
}
\caption{An illustration of the fact that the two types of geometric hyperplanes of $S_{(2)}$ ({\it right}) can be regarded as blow-ups of the Veldkamp lines of $S_{(1)}$ ({\it left}), or, {\it vice versa}, that a projection of a type two ({\it top right}) or type one ({\it bottom right}) hyperplane of $S_{(2)}$ onto a line of $S_{(2)}$ can be viewed, respectively, as the ordinary ({\it top left}) or an extraordinary ({\it bottom left}) Veldkamp line of $S_{(1)}$.}
\end{figure}
 
In order to get a deeper insight into the nature of the two (ordinary) types of geometric hyperplanes of $S_{(2)}$, we shall have a look at Veldkamp lines of $S_{(1)}$. These are also of two different kinds, which we shall call ordinary (abbreviated as `ord' in Table 1) and extraordinary (`ext') ones.
The Veldkamp line of ordinary type is just a single one, namely that comprising all the three points of $S_{(1)}$ --- as portrayed in Figure 2, top left, where the points are distinguished by different color. A Veldkamp line of extraordinary type features the extraordinary hyperplane (that is the whole $S_{(1)}$) and a point counted twice --- see Figure 2, bottom left; clearly, there are three of them.
Now, $S_{(2)}$ is endowed with two different {\it spreads of lines}, i.\,e. sets of three pairwise disjoint lines that partition the point-set of $S_{(2)}$; in our diagrammatical representation of $S_{(2)}$, one spread consists of three `horizontal' lines and the other comprises three `vertical' lines. Let us take a geometric hyperplane of type two and project its points along the three lines of one spread onto a line of the second spread. If we keep the distinction between the three lines of the second spread, what we get is nothing but the ordinary Veldkamp line of $S_{(1)}$. This fact is portrayed in Figure 2, top, where the projection is made along the lines of the `vertical' spread. Making the same projection with a type one hyperplane results in an extraordinary Veldkamp line of $S_{(1)}$ --- see Figure 2, bottom part. We further see that {\it six} different $S_{(2)}$-hyperplanes of type two are projected into one and the same {\it ordinary} Veldkamp line of $S_{(1)}$, whereas it is {\it three} 
distinct $S_{(2)}$-hyperplanes of type one that map into the same {\it extraordinary} Veldkamp line of $S_{(1)}$. Reversing our reasoning, by blowing up $S_{(1)}$ to $S_{(2)}$ each of the three extraordinary Veldkamp lines of $S_{(1)}$ generates three geometric hyperplanes of $S_{(2)}$ of type one, whilst the single ordinary one gives rise to all the six hyperplanes of type two. We have thus explicitly demonstrated the fact that the types of geometric hyperplanes, that is Veldkamp {\it points}, of the Segre variety $S_{(2)}$, as well as their cardinalities, can be fully derived/recovered from the properties of Veldkamp {\it lines} of the Segre variety of the preceding rank, $S_{(1)}$. One of the main foci of this paper is to justify a conjecture that this remarkable relation between Veldkamp points of $S_{(N)}$ and Veldkamp lines of $S_{(N-1)}$ holds, in fact, for any positive integer $N \geq 2$ by illustrating its properties in sufficient detail on the two subsequent cases, viz. $N=3$ and $N=4$.

\begin{table}[b]
\begin{center}
\caption{The types of ordinary Veldkamp lines of $S_{(2)}$. The first column gives the type, the next two columns tell us about how many points and lines belong to all the three geometric hyperplanes a line of the given type consists of, then we learn about the line's composition and, finally, the last column lists cardinalities for each type.} 
\vspace*{0.2cm}
{\begin{tabular}{|c|c|c|c|c|c|} \hline \hline
\multicolumn{1}{|c|}{}  &  \multicolumn{2}{|c|}{} & \multicolumn{2}{|c|}{}  &  \multicolumn{1}{|c|}{}\\
\multicolumn{1}{|c|}{} &   \multicolumn{2}{|c|}{Core} & \multicolumn{2}{|c|}{Comp'n}  &\multicolumn{1}{|c|}{} \\
 \cline{2-5}
Tp & Ps & Ls  & $H_1$ & $H_2$ & Crd  \\
\hline
1  &  3   & 1     & 3     & --    &  6 \\
\hline
2  &  2   & 0     & 2     & 1     &  18 \\
\hline
3  &  1   & 0     & 1     & 2     &  9  \\
\hline
4  &  0   & 0     & --    & 3     &  2 \\
 \hline \hline
\end{tabular}}
\end{center}
\end{table}

It is now obvious that in order to tackle properly the next ($N=3$) case in the hierarchy, we have to ascertain how many types of Veldkamp lines of $S_{(2)}$ we have and what their properties are. Using our descriptive representation of geometric hyperplanes of $S_{(2)}$ (Figure 1), this is quite a straightforward task. For ordinary Veldkamp lines the corresponding information is collected in Table 2 and, in a visual form, in Figure 3, top. There are four types of them and their total number is 35; type one and type four are so-called homogenous Veldkamp lines as they both feature hyperplanes of the same type. As for extraordinary Veldkamp lines, we find only two types, namely the extraordinary hyperplane and $H_1$ counted twice (type I) and the extraordinary hyperplane and $H_2$ counted twice (type II), totaling to 15 -- as depicted in the upper part of Figure 4, left and right, respectively. 

We shall finalize this section by observing that the ordinary part of Veldkamp space of $S_{(2)}$, ord-$\mathcal{V}(S_{(2)})$, features 15 points (ordinary geometric hyperplanes of $S_{(2)}$) and 35 lines (ordinary Veldkamp lines of $S_{(2)}$), being isomorphic to $PG(3,2)$.
The nine hyperplanes of type one then correspond to the nine points lying on a hyperbolic quadric $\mathcal{Q}_0^{+}(3,2) \subset PG(3,2)$.

\begin{figure}[t]
\centerline{\includegraphics[width=11truecm,clip=]{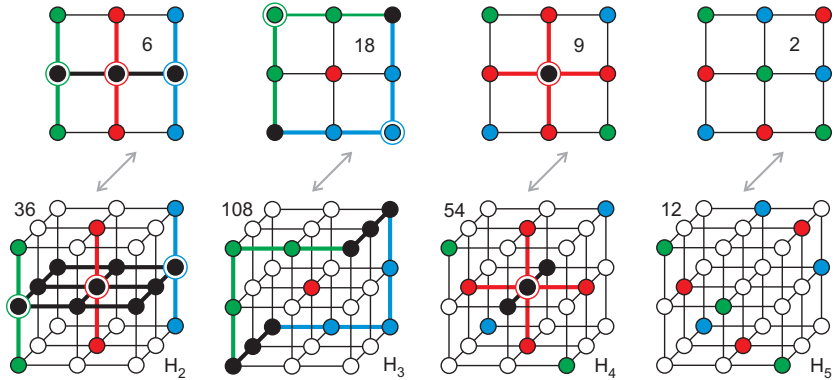}}
\caption{{\it Top:} -- A descriptive illustration of the structure of the four distinct types  (1 to 4, left to right) of ordinary Veldkamp lines of $S_{(2)}$ (Table 2); the three geometric hyperplanes comprising a Veldkamp line are distinguished by different colors, with the points and lines shared by all of them being colored black. {\it Bottom:} -- The four distinct types of geometric hyperplanes of $S_{(3)}$, as well as the number of copies per each type, we get by blowing-up Veldkamp lines of $S_{(2)}$ of the type shown above the particular subfigure.}
\end{figure}

\begin{figure}[pth!]
\centerline{\includegraphics[width=5.3truecm,clip=]{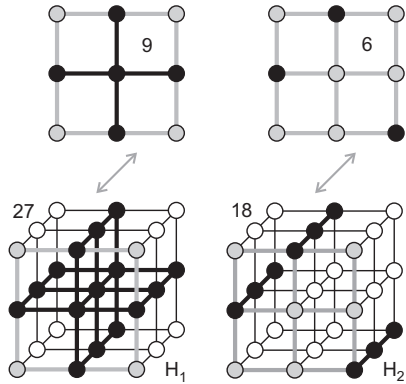}}
\caption{The same as in Figure 3, but for extraordinary Veldkamp lines of $S_{(2)}$ ({\it top}) and their $S_{(3)}$ blown-up cousins ({\it bottom}).}
\end{figure}

\subsection{\textbf{\textit{S}}$_{(3)}$ and its Veldkamp space}\label{S3} 

The Segre variety $S_{(3)}$, known also as the smallest slim dense near hexagon \cite{bru} or the Gray configuration \cite{mpw}, is a point-line incidence structure having 27 points and the same number of lines, with three points per line and the same number of lines through a point. It possesses nine $S_{(2)}$'s, arranged into three distinct triples of pairwise disjoint members, each of which partitions the point-set. Analogously,  $S_{(3)}$ contains three {\it distinguished} spreads of lines; a distinguished spread of lines is a set of nine mutually skew lines such that each line is incident with all the three $S_{(2)}$'s from a given triple.\footnote{In analogy to $S_{(2)}$, and to be compatible with the approach advanced in \cite{grsa}, we shall employ a handy descriptive representation of $S_{(3)}$ as a $3 \times 3 \times 3$ grid (see Figure 1 of \cite{grsa}).} As already mentioned, $S_{(3)}$, as an abstract point-line incidence structure, was thoroughly analyzed in \cite{grsa} where it was found that ord-$\mathcal{V}(S_{(3)}) \cong PG(7,2)$. As $PG(7,2)$ has 255 points and 10\,795 lines (see, e.\,g., \cite{ht}), $S_{(3)}
$ features a total of 255 (ordinary) geometric hyperplanes and 10\,795 (ordinary) Veldkamp lines.

\begin{figure}[t]
\centerline{\includegraphics[width=12.0truecm,clip=]{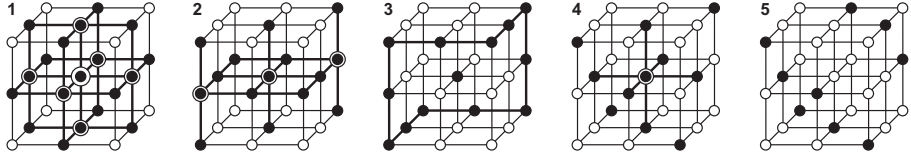}}
\caption{A diagrammatic $3 \times 3 \times 3$ grid representation of the types of (ordinary) geometric hyperplanes of $S_{(3)}$ (after \cite{grsa}, Figure 2). As in Figure 1, deep points are encircled. Note that one of the seven deep points of a type-one hyperplane stands on a different footing than the others, as each line passing through it features deep points; we shall call it the deepest point.}
\end{figure}

\begin{table}[pth!]
\begin{center}
\caption{The 5 types of (ordinary) geometric hyperplanes of the Segre variety $S_{(3)}$. As in Table 1, the first column gives the type (`Tp') of a hyperplane, which is followed by the number of points (`Ps') and lines (`Ls') it contains, and the number of points of given order. The next three columns tell us about how many of 9 $S_{(2)}$'s are fully located (`D') in the hyperplane and/or share with it a hyperplane of type $H_1$  or $H_2$ (see Table 1). The VL-column lists the types of (ordinary and extraordinary) Veldkamp lines of $S_{(2)}$  we get by projecting a hyperplane of the given type into an $S_{(2)}$ along the lines of all three distinguished spreads. Finally, for each hyperplane type we give its cardinality (`Crd'), the corresponding large orbit  of $2 \times 2 \times 2$ arrays over $GF(2)$ (`BS') taken from Table 3 of \cite{brst}, and its weight, or rank in the language of \cite{brst} (`W').} 
\vspace*{0.3cm}
{\begin{tabular}{|c|c|c|c|c|c|c|c|c|c|c|c|c|c|} \hline \hline
\multicolumn{1}{|c|}{} & \multicolumn{1}{|c|}{} & \multicolumn{1}{|c|}{}  &  \multicolumn{4}{|c|}{} & \multicolumn{3}{|c|}{}  &\multicolumn{1}{|c|}{} &\multicolumn{1}{|c|}{} &\multicolumn{1}{|c|}{}  & \multicolumn{1}{|c|}{}\\
\multicolumn{1}{|c|}{} & \multicolumn{1}{|c|}{} &
\multicolumn{1}{|c|}{}  &  \multicolumn{4}{|c|}{Points of
Order} & \multicolumn{3}{|c|}{$S_{(2)}$'s of Type} & \multicolumn{1}{|c|}{} & \multicolumn{1}{|c|}{} &\multicolumn{1}{|c|}{} 
& \multicolumn{1}{|c|}{}\\
 \cline{4-10}
Tp & Ps & Ls  & ~0~ & ~1~ & 2 & ~3~ & ~D~ & $H_1$ & $H_2$ & VL & Crd & BS & W \\
\hline
1 & 19 & 15  & 0 & 0 & 12   & 7 & 3 & 6 & 0 & I      & 27  & 2  & 1 \\
\hline
2 & 15 & 9   & 0 & 6 &  6   & 3 & 1 & 6 & 2 & II,\,1 & 54  & 3  & 2 \\
\hline
3 & 13 & 6   & 1 & 6 & 6    & 0 & 0 & 6 & 3 &  2     & 108 & 4  & 2 \\
\hline
4 & 11 & 3   & 4 & 6 & 0    & 1 & 0 & 3 & 6 &  3     & 54  & 5  & 3 \\
\hline
5 & 9  & 0   & 9 & 0 & 0    & 0 & 0 & 0 & 9 &  4     & 12  & 6  & 3 \\
\hline \hline
\end{tabular}}
\end{center}
\end{table}

Let us first focus briefly on geometric hyperplanes of $S_{(3)}$, referring the interested reader to \cite{grsa} (as well as to Sec.\,6 of \cite{hos2}) for more details. Disregarding the extraordinary one, there are five distinct types of 
them, their basic properties being listed in Table 3 and their representatives depicted in Figure 5; we note that a type one hyperplane is the singular one, whereas that of type five is an ovoid. Type one hyperplanes, as already pointed out in \cite{grsa}, play a distinguished role in the sense that a hyperplane of any other type can be expressed as the Veldkamp sum of several type-one hyperplanes; the smallest such number is then called the {\it weight} of the hyperplane.
A type two hyperplane also differs from the rest in being a blow-up of two different kinds of Veldkamp lines of $S_{(2)}$, viz. those of type 1 ordinary (Figure 3, leftmost) and those of type II extraordinary (Figure 4, right), with its 36/18 copies originating from the former/latter process; hyperplanes of any of the remaining four types originate --- as it can easily be verified from the corresponding parts of Figure 3 and Figure 4  --- from $S_{(2)}$-Veldkamp-lines of the same type.

\begin{table}[pth!]
\begin{center}
\caption{The types of ordinary Veldkamp lines of $S_{(3)}$, in the notation set up in Table 2. A slightly modified reproduction of Table 2 from \cite{grsa}, where the reader is referred to for detailed explanation of the meaning of all subscripts associated to some cardinalities of core's points and lines (columns two and three).} \vspace*{0.6cm}
{\begin{tabular}{|r|l|l|c|c|c|c|c|r|} \hline \hline
\multicolumn{1}{|c|}{}  &  \multicolumn{2}{|c|}{} & \multicolumn{5}{|c|}{}  &  \multicolumn{1}{|c|}{}\\
\multicolumn{1}{|c|}{} &   \multicolumn{2}{|c|}{Core} & \multicolumn{5}{|c|}{Composition}  &\multicolumn{1}{|c|}{} \\
 \cline{2-8}
Tp & ~Ps & ~Ls  & $H_1$ & $H_2$ & $H_3$ & $H_4$ & $H_5$ & Crd  \\
\hline
1 &  ~15 & ~11  & 3 & -- & -- & -- & -- &  27 \\
\hline
2 & ~13 & ~8  & 2 & 1 & -- & -- & -- &  162 \\
\hline
3 & ~12 & ~6  & 2 & -- & 1 & -- & -- & 108  \\
\hline
4 &  ~11 & ~7  & 1 & 2 & -- & -- & -- &  81 \\
\hline
5 &  ~10 & ~4  & 1 & 1 & 1 & -- & -- & 648  \\
\hline
6 &  ~9 & ~6  & -- & 3 & -- & -- & -- & 18  \\
7 &  ~9 & ~4  & 1 & -- & 2 & -- & -- &  324 \\
8 &  ~$9_{(2)}$ & ~3c  & 1 & 1 & -- & 1 & -- &  324 \\
9 &  ~9 & ~3  & 1 & -- & 2 & -- & -- &  324 \\
10 &  ~9 & ~3p  & -- & 3 & -- & -- & -- &  18 \\
11 &  ~$9_{(3)}$ & ~3c  & -- & 3 & -- & -- & -- &  108 \\
\hline
12 &  ~8 & ~3  & -- & 2 & 1 & -- & -- & 648  \\
13 &  ~8 & ~2  & 1 & -- & 1 & 1 & -- &  648 \\
\hline
14 &  ~7 & ~3  & 1 & -- & -- & 2 & -- &  27 \\
15 &  ~7 & ~2p  & -- & 1 & 2 & -- & -- & 162  \\
16 & ~$7_{(2)}$  & ~2c   & -- & 1 & 2 & -- & -- & 324  \\
17 & ~$7_{(3)}$  & ~2c   & -- & 1 & 2 & -- & -- & 324  \\
18 &  ~$7_{[2]}$ & ~1  & -- & 2 & -- & 1 & -- & 162  \\
19 &  ~$7_{[1]}$ & ~1  & -- & 1 & 2 & -- & -- & 324  \\
20 &  ~7 & ~0  & 1 & -- & 1 & -- & 1 & 108  \\
21 &  ~7 & ~0  & 1 & -- & -- & 2 & -- &  108 \\
\hline
22 &  ~6 & ~2c  & -- & 1 & 1 & 1 & -- & 648  \\
23 &  ~6 & ~2p  & -- & -- & 3 & -- & --  & 108  \\
24 &  ~6 & ~1  & -- & -- & 3 & -- & -- &  648 \\
25&  ~$6_{[3]}$ & ~0  & 1 & -- & -- & 1 & 1 & 216  \\
26 &  ~$6_{[2]}$ & ~0  & -- & 2 & -- & -- & 1 & 108  \\
27 &  ~$6_{[1]}$ & ~0  & -- & 1 & 1 & 1 & -- & 648  \\
28 &  ~$6_{[0]}$ & ~0  & -- & -- & 3 & -- & -- & 36  \\
\hline
29 &  ~$5_{[1]}$ & ~1  & -- & 1 & -- & 2 & -- &  162 \\
30 &  ~$5_{[0]}$ & ~1  & -- & -- & 2 & 1 & -- & 648  \\
31 &  ~$5_{(2)}$ & ~0  & -- & 1 & 1 & -- & 1 &  324 \\
32 &  ~$5_{(1)}$ & ~0  & -- & 1 & -- & 2 & -- &  324 \\
33 &  ~$5_{(0)}$ & ~0  & -- & -- & 2 & 1 & -- & 648  \\
\hline
34 &  ~4 & ~0  & -- & -- & 2 & -- & 1 & 324  \\
35 &  ~$4_{(3:1)}$ & ~0  & -- & -- & 1 & 2 & -- & 216  \\
36 &  ~$4_{(2:2)}$ & ~0  & -- & -- & 1 & 2 & -- & 324  \\
\hline
37 &  ~3 & ~1  & -- & -- & -- & 3 & -- & 54  \\
38 &  ~$3_{[1]}$ & ~0  & -- & 1 & -- & -- & 2 &  54 \\
39 &  ~$3_{[0]}$ & ~0  & -- & -- & 1 & 1 & 1 & 216  \\
\hline
40 &  ~2 & ~0  & -- & -- & -- & 2 & 1 &  108 \\
\hline
41 &  ~0 & ~0  & -- & -- & -- & -- & 3 & 4  \\
 \hline \hline
\end{tabular}}
\end{center}
\end{table}

Regarding hyperplanes as points of the $PG(7,2)$ ($\cong$ ord-$\mathcal{V}(S_{(3)})$), it can be demonstrated that those points that correspond to hyperplanes of type one, two and four, and whose number totals to 135, all lie on a certain $\mathcal{Q}_0^{+}(7,2) \subset PG(7,2)$; note that these are exactly the three types whose members feature points of maximum order. It is also worth mentioning here that these three hyperplane types in their totality correspond to the image, furnished by the Lagrangian Grassmannian of $LG(3,6)$ type, of the set of 135 maximal subspaces of $\mathcal{W}(5,2)$ \cite{hsl,leps}.

Next, let us have an informative look at the (ordinary) Veldkamp lines of $S_{(3)}$. As already mentioned, there are 10\,795 of them and, as found in \cite{grsa}, they fall into 41 distinct types. The basic combinatorial characteristics of all the types are listed in Table 4 and a representative of each type is provided by Figure 6.
Concerning extraordinary Veldkamp lines of $S_{(3)}$, there are obviously as many types of them as there are types of ordinary geometric hyperplanes of $S_{(3)}$; we depict their representatives in Figure 7.

\begin{figure}[pth!]
\centerline{\includegraphics[width=5.5truecm,clip=]{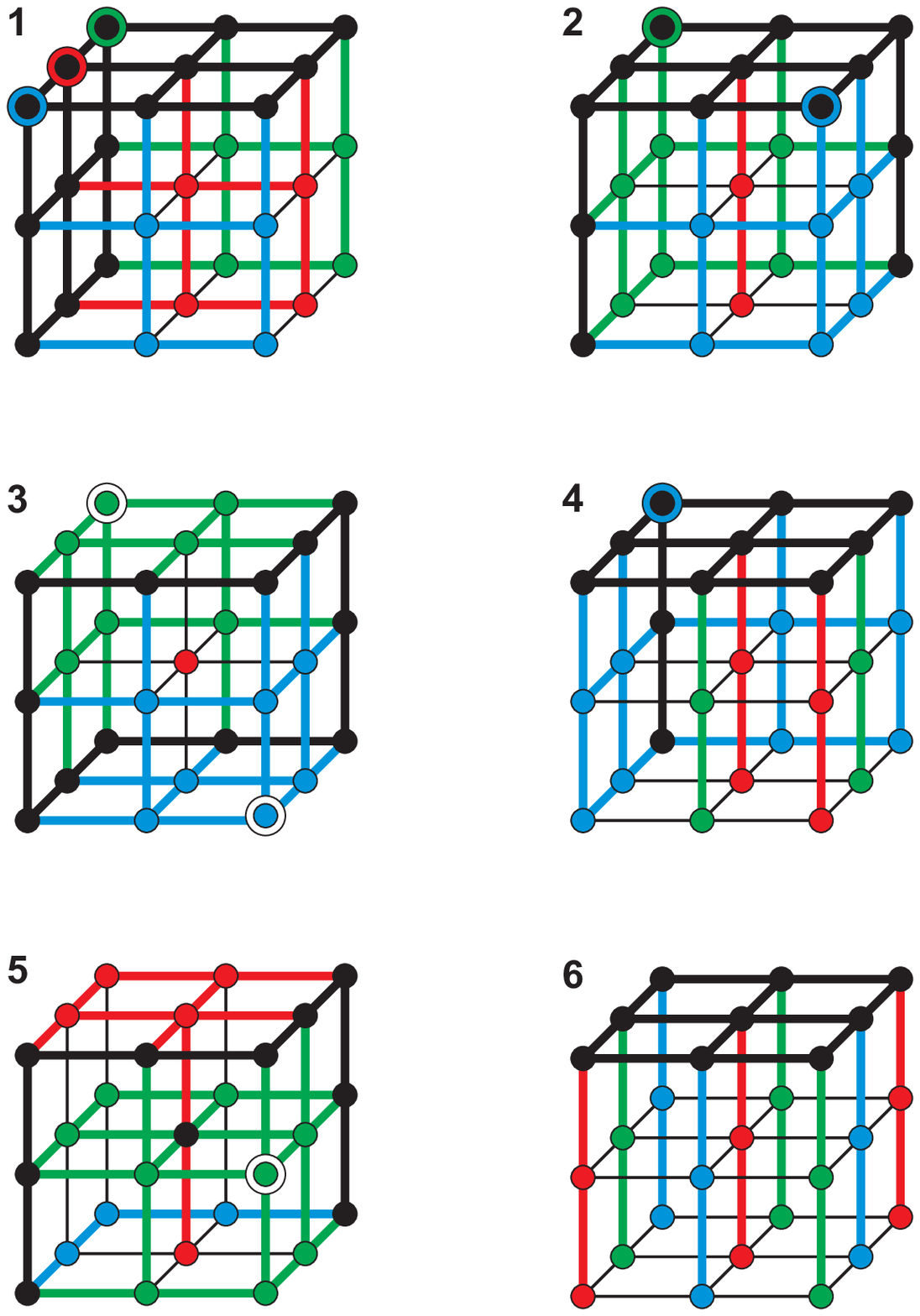}\hspace*{1.7cm}\includegraphics[width=5.5truecm,clip=]{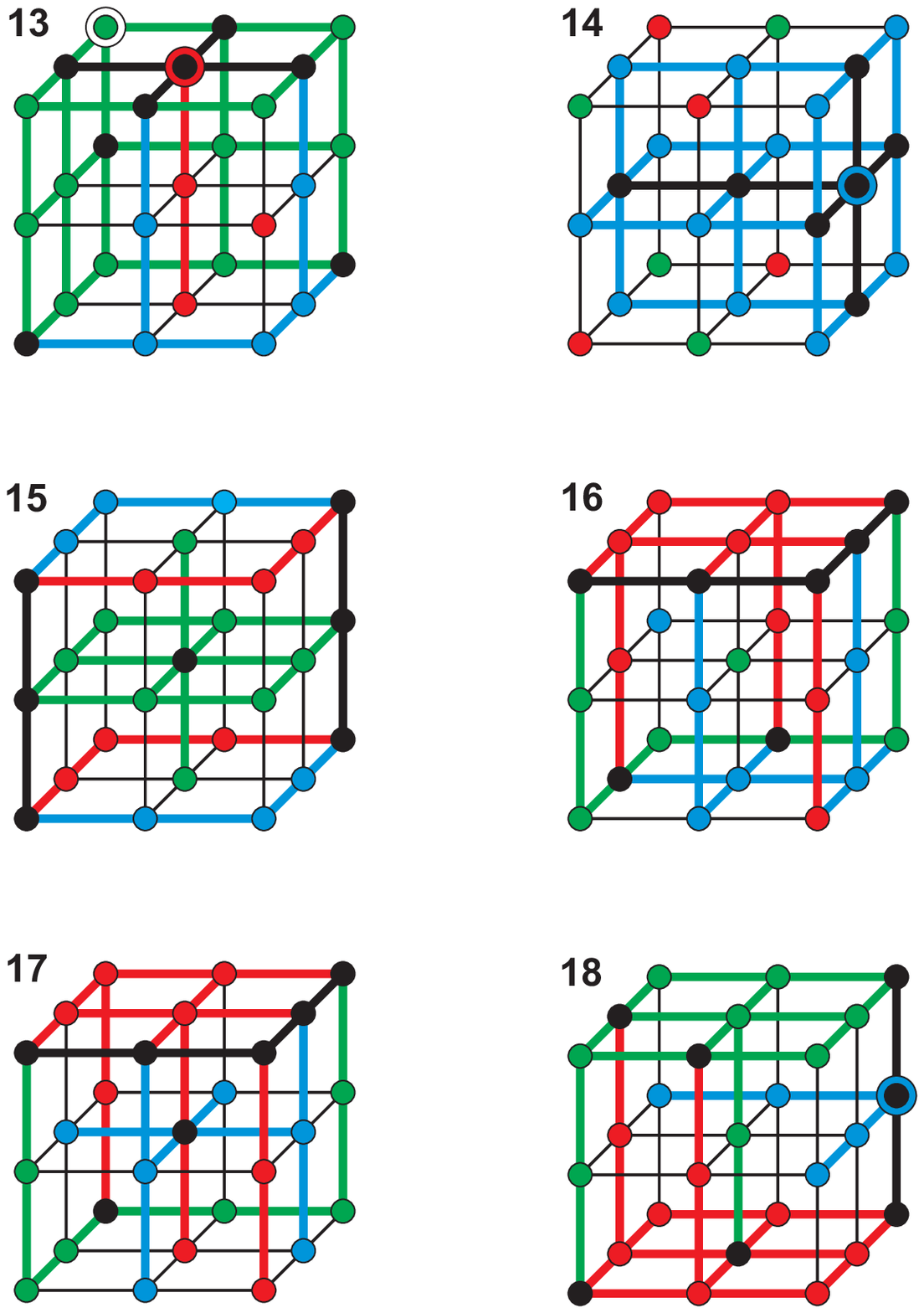}}
\vspace*{0.7cm}
\centerline{\includegraphics[width=5.5truecm,clip=]{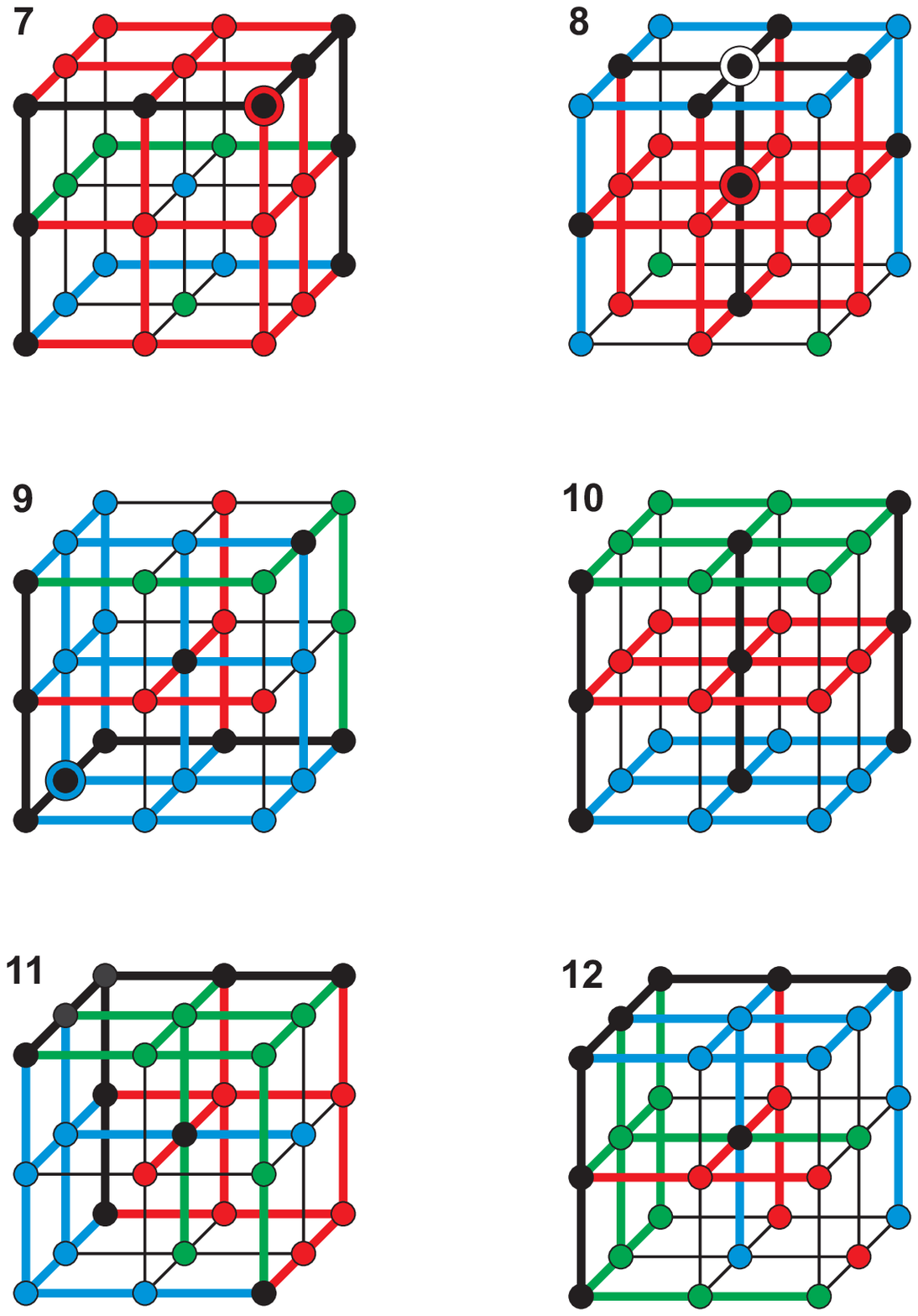}\hspace*{1.7cm}\includegraphics[width=5.5truecm,clip=]{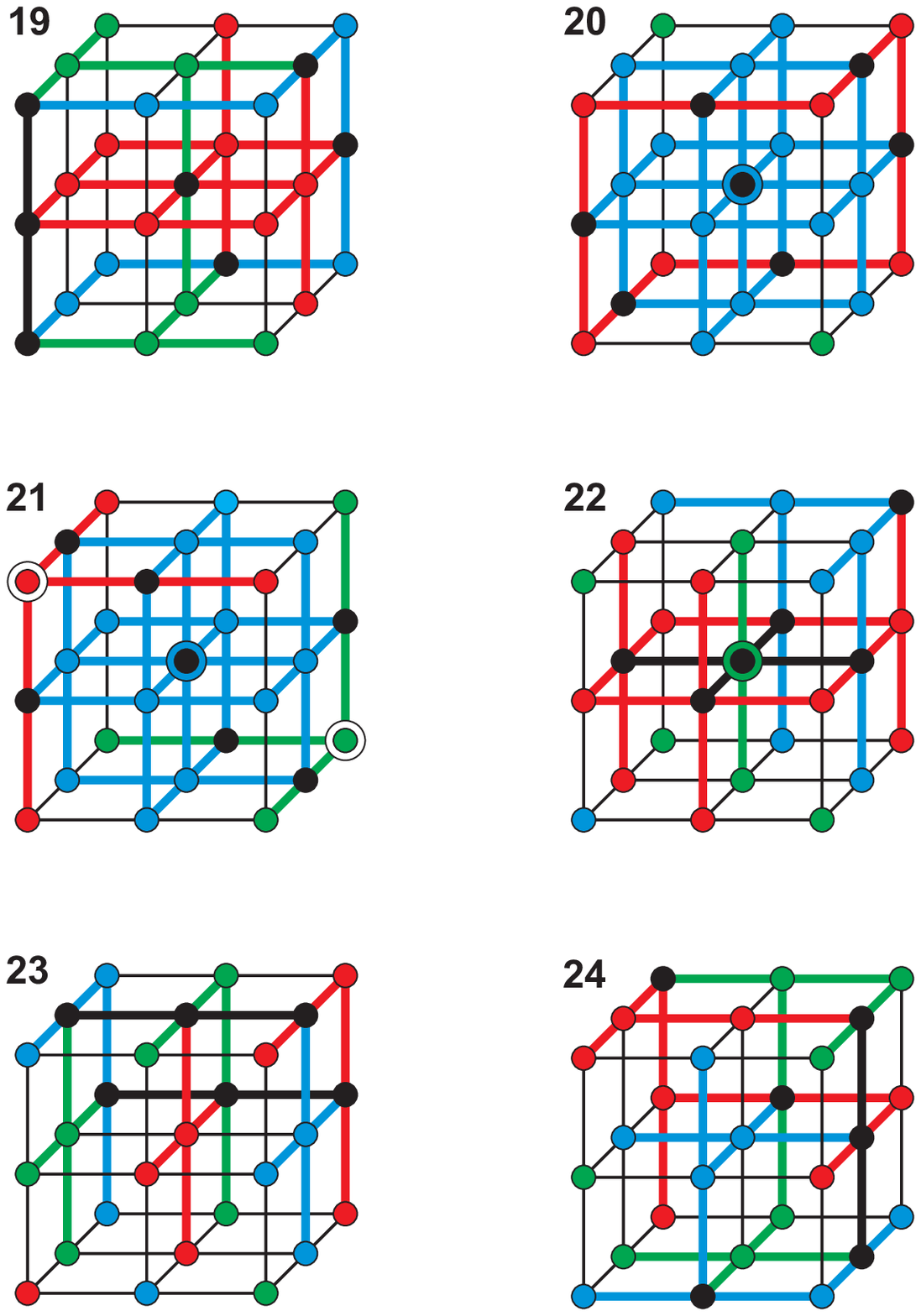}}
\caption{A diagrammatic illustration, on a $3 \times 3 \times 3$-grid, of representatives of all 41 types of ordinary Veldkamp lines of $S_{(3)}$. As in Figure 3, top, the three geometric hyperplanes forming a Veldkamp line are distinguished by different coloring, with the points and lines shared by all of them (the core) being colored black.}
\end{figure}
\clearpage
\addtocounter{figure}{-1}
\begin{figure}[pth!]
\centerline{\includegraphics[width=5.5truecm,clip=]{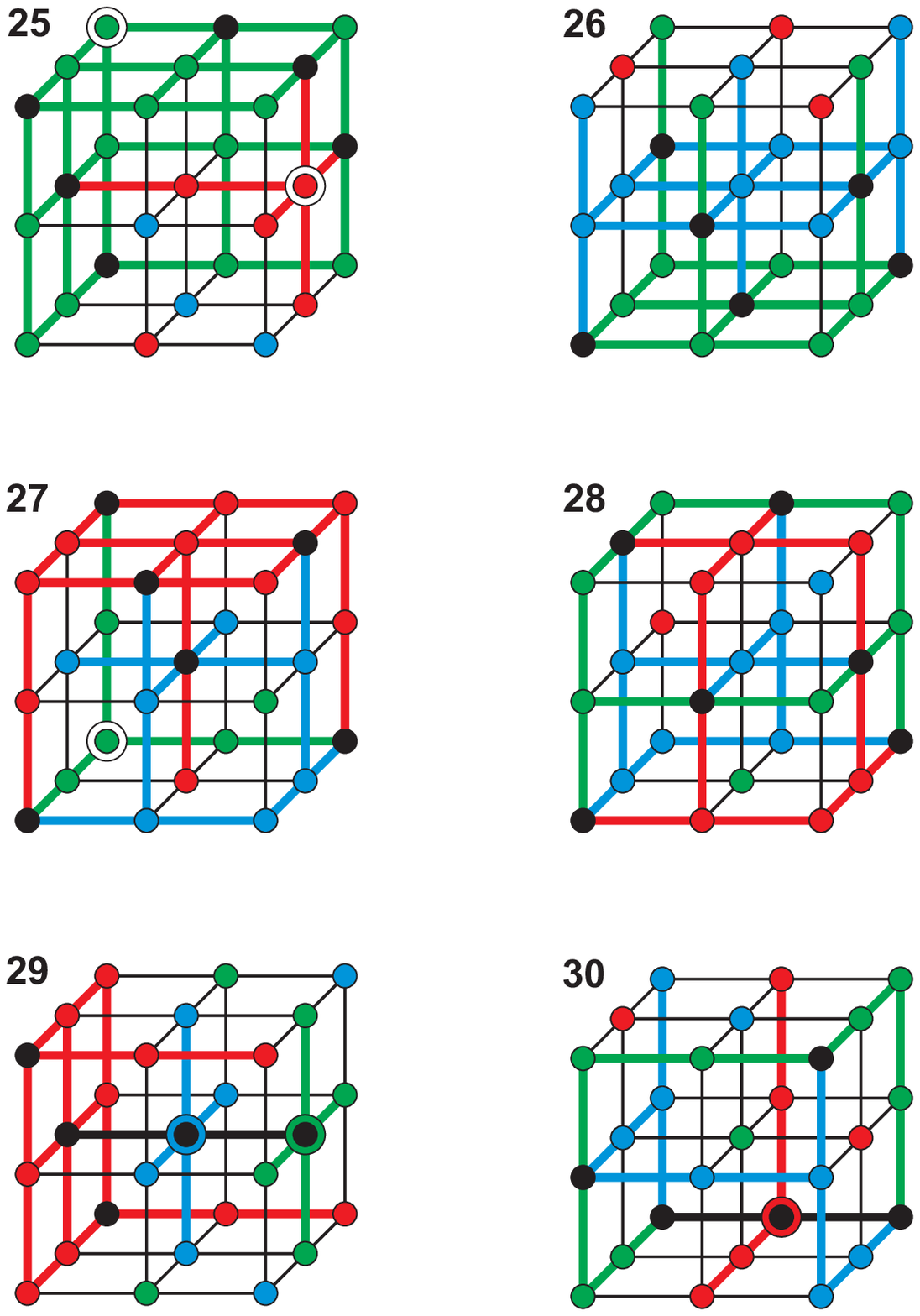}\hspace*{1.7cm}\includegraphics[width=5.5truecm,clip=]{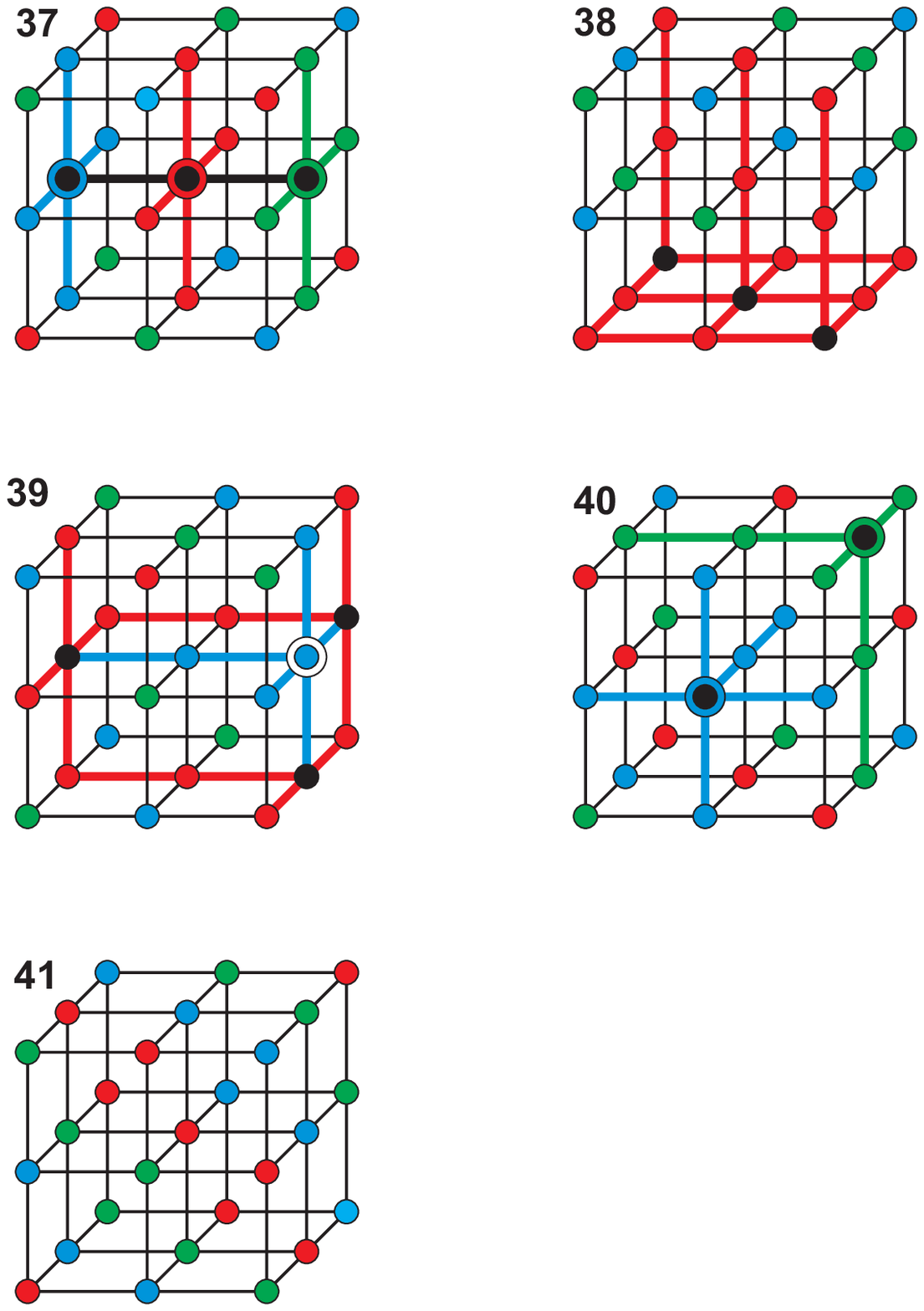}}
\vspace*{0.7cm}
\centerline{\includegraphics[width=5.5truecm,clip=]{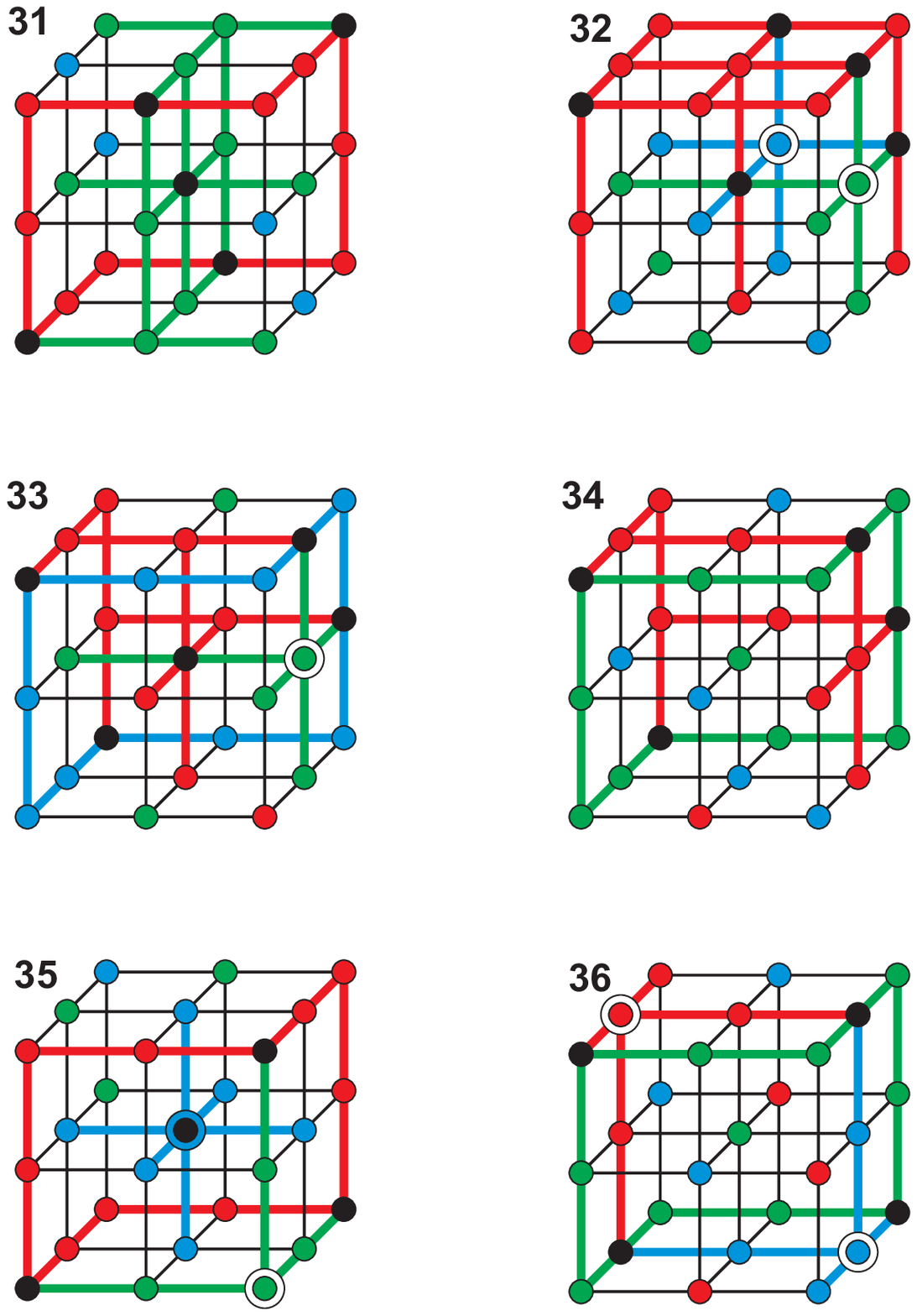}\hspace*{7.1cm}}
\caption{(Continued.)}
\end{figure}

\begin{figure}[pth!]
\centerline{\includegraphics[width=13truecm,clip=]{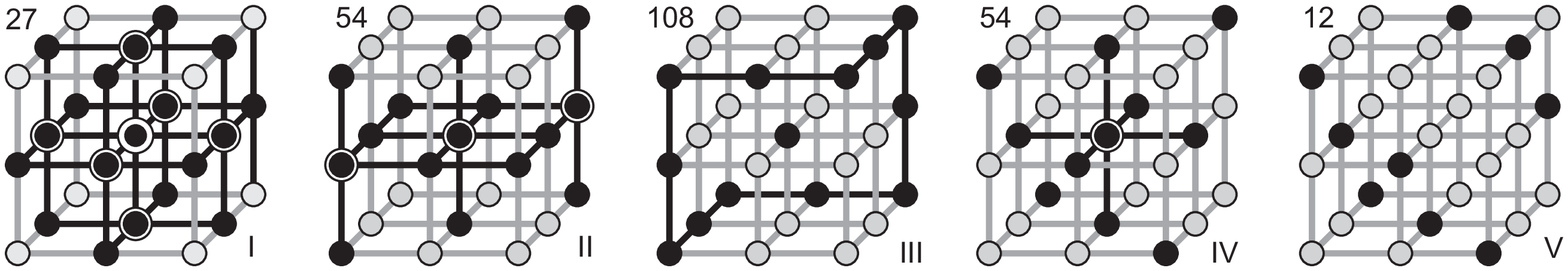}}
\caption{A diagrammatic $3 \times 3 \times 3$ grid representation of the five types of extraordinary Veldkamp lines of $S_{(3)}$. Roman numerals stand for the type, whilst Arabic ones denote the number of members for each type.}
\end{figure}

\subsection{\textbf{\textit{S}}$_{(4)}$ and its geometric hyperplanes}\label{S4} 
The Segre variety $S_{(4)}$, known also as the smallest slim dense near octagon \cite{bru}, is a point-line incidence structure having 81 points and 108 lines, with three points per line and four lines through a point. It possesses 12 $S_{(3)}$'s, arranged into four distinct triples of pairwise disjoint members, each of which partitions the point-set. Analogously,  $S_{(4)}$ 
contains four {\it distinguished} spreads of lines; a distinguished spread of lines is a set of 27 mutually skew lines such that each line is incident with all the three $S_{(3)}$'s from a given triple.

\begin{table}[pht!]
\begin{center}
\caption{The 29 types of  geometric hyperplanes of $S_{(4)}$; also shown is partition of hyperplane types into 15 classes according to the number of points/lines. As in corresponding Tables 1 and 3, one first gives the type (`Tp') of a hyperplane, then the number of points (`Ps') and lines (`Ls') it contains, and the number of points of given order. The next six columns tell us about how many of 12 $S_{(3)}$'s are fully located (`D') in the hyperplane and/or share with it a hyperplane of type $H_i$ (see Table 3 above). The VL-column lists the types of ordinary and/or extraordinary Veldkamp lines of $S_{(3)}$  we get by projecting a hyperplane of the given type into $S_{(3)}$'s along the lines of all four distinguished spreads. Finally, for each hyperplane type we give its cardinality (`Crd'), the corresponding large orbit  of $2 \times 2 \times 2 \times 2$ arrays over $GF(2)$ (`BS') taken from Table 5 of \cite{brst}, and its weight/rank (`W'). }
{\begin{tabular}{|r|c|c|c|c|c|c|c|c|c|c|c|c|c|r|r|r|c|} \hline \hline
\multicolumn{1}{|c|}{} & \multicolumn{1}{|c|}{} & \multicolumn{1}{|c|}{}  &  \multicolumn{5}{|c|}{}                        & \multicolumn{6}{|c|}{}                      & \multicolumn{1}{|c|}{}  & \multicolumn{1}{|c|}{} & \multicolumn{1}{|c|}{} & \multicolumn{1}{|c|}{} \\
\multicolumn{1}{|c|}{} & \multicolumn{1}{|c|}{} & \multicolumn{1}{|c|}{}  &  \multicolumn{5}{|c|}{$\#$ of Points of Order} & \multicolumn{6}{|c|}{$\#$ of $S_{(3)}$'s of Type} & \multicolumn{1}{|c|}{}  & \multicolumn{1}{|c|}{} & \multicolumn{1}{|c|}{} & \multicolumn{1}{|c|}{} \\
 \cline{4-14}
Tp & Ps & Ls  & 0 & 1 & 2 & 3 & 4 &  D  & $H_1$ & $H_2$ & $H_3$ & $H_4$ &  $H_5$ & VL         & Crd  & BS & W\\
\hline
1  &  65 & 76   &  0  &  0  &  0  &  32 & 33  &  4   &    8   &   0   &   0   &   0   &   0    & I       & 81   & 2  & 1 \\ 
\hline
2  &  57 & 60   &  0  &  0  & 12  &  24 & 21  &  2   &    6   &   4   &   0   &   0   &   0    &   II,\,1        & 324  & 3  & 2 \\
\hline
3  &  53 & 52   &  0  &  2  & 12  &  26 & 13  &  1   &    6   &   3   &   2   &   0   &   0    &   III,\,2        & 1296 & 4  & 2 \\
\hline 
4  &  51 & 48   &  1  &  0  &  12 &  32 &  6  &  0   &    8   &   0   &   4   &   0   &   0    &   3        & 648  & 5  & 2 \\
\hline
5  &  49 & 44   &  0  &  8  & 12  &  16 & 13  &  1   &    3   &   6   &   0   &   2   &   0    &   IV,\,4        & 648  & 6  & 3 \\
\hline
6  &  47 & 40   &  0  &  4  &  18 &  20 &  5  &  0   &    4   &   4   &   4   &   0   &   0    &   5        & 3888 & 11 & 3 \\
\hline
7  &  45 & 36   &  0  & 18  &  0  &  18 &  9  &  1   &    0   &   9   &   0   &   0   &   2    &   V,\,6        & 144  & 7  & 3 \\
8  &  45 & 36   &  0  &  0  &  36 &  0  &  9  &  0   &    0   &   12  &   0   &   0   &   0    &  10        & 108  & 17 & 4 \\
9  &  45 & 36   &  2  &  4  &  18 & 16  &  5  &  0   &    4   &   2   &   4   &   2   &   0    &   7,\,8    & 3888 & 8  & 3 \\
10 &  45 & 36   &  0  &  6  &  18 & 18  &  3  &  0   &    3   &   3   &   6   &   0   &   0    &   9,\,11   & 2592 & 9  & 3 \\
\hline
11 &  43 & 32   &  1  &  8  &  18 & 12  &  4  &  0   &    2   &   4   &   4   &   2   &   0    &   12,\,13  & 7776 & 12  & 3 \\ 
\hline
12 &  41 & 28   &  8  &  0  &  24 &  0  &  9  &  0   &    4   &   0   &   0   &   8   &   0    &   14       & 162  & 14 & 4 \\ 
13 &  41 & 28   &  0  &  12 &  18 &  8  &  3  &  0   &    0   &   6   &   4   &   2   &   0    &   15,\,18  & 1944 & 19   & 4 \\ 
14 &  41 & 28   &  0  &  14 &  12 & 14  &  1  &  0   &    1   &   3   &   6   &   2   &   0    &   17,\,21  & 2592 & 15 & 4 \\ 
15 &  41 & 28   &  2  &  8  &  18 & 12  &  1  &  0   &    1   &   3   &   7   &   0   &   1    &   16,\,20  & 2592 & 10   & 3 \\ 
16 &  41 & 28   &  0  &  8  &  24 &  8  &  1  &  0   &    0   &   4   &   8   &   0   &   0    &   19       & 1944 & 20  & 4 \\ 
\hline
17 &  39 & 24   &  4  &  12 & 12  &  8  &  3  &  0   &    1   &   3   &   3   &   4   &   1    &   22,\,25  & 5184 & 16 & 4 \\ 
18 &  39 & 24   &  3  &  12 & 12  & 12  &  0  &  0   &    0   &   4   &   6   &   0   &   2    &   23,\,26  & 1296 & 22 & 4 \\ 
19 &  39 & 24   &  1  &  12 & 18  &  8  &  0  &  0   &    0   &   2   &   8   &   2   &   0    &   24,\,27  & 7776 & 23 & 4 \\ 
20 &  39 & 24   &  3  &  0  & 36  &  0  &  0  &  0   &    0   &   0   &  12   &   0   &   0    &   28       & 216  & 13 & 3 \\ 
\hline
21 &  37 & 20   &  4  & 16  & 12  &  0  &  5  &  0   &    0   &   4   &   0   &   8   &   0    &   29       & 972  & 18 & 4 \\ 
22 &  37 & 20   &  4  & 14  & 12  &  6  &  1  &  0   &    0   &   2   &   5   &   4   &   1    &   30--32   & 7776 & 21 & 4 \\ 
23 &  37 & 20   &  3  & 12  & 18  &  4  &  0  &  0   &    0   &   0   &   8   &   4   &   0    &   33       & 3888 & 24 & 4 \\ 
\hline
24 &  35 & 16   &  4  & 20  & 6   &  4  &  1  &  0   &    0   &   0   &   4   &   8   &   0    &   35       & 1296 & 28 & 5 \\ 
25 &  35 & 16   &  7  & 12  & 12  &  4  &  0  &  0   &    0   &   0   &   6   &   4   &   2    &   34,\,36  & 3888 & 25 & 4 \\ 
\hline
26 &  33 & 12   &  12 &  12 & 6   &  0  &  3  &  0   &    0   &   2   &   0   &   6   &   4    &   37,\,38  & 648  & 26 & 5 \\
27 &  33 & 12   &  11 &  12 & 6   &  4  &  0  &  0   &    0   &   0   &   4   &   4   &   4    &   39       & 1296 & 29 & 5 \\
\hline
28 &  31 &  8   &  13 &  16 & 0   &  0  &  2  &  0   &    0   &   0   &   0   &   8   &   4    &   40       & 648  & 27 & 5 \\
\hline
29 &  27 &  0   &  27 &  0  & 0   &  0  &  0  &  0   &    0   &   0   &   0   &   0   &   12   &   41       & 24   & 30 & 6 \\
\hline \hline
\end{tabular}}
\end{center}
\end{table}

\begin{figure}[t]
\centerline{\includegraphics[width=3.8truecm,clip=]{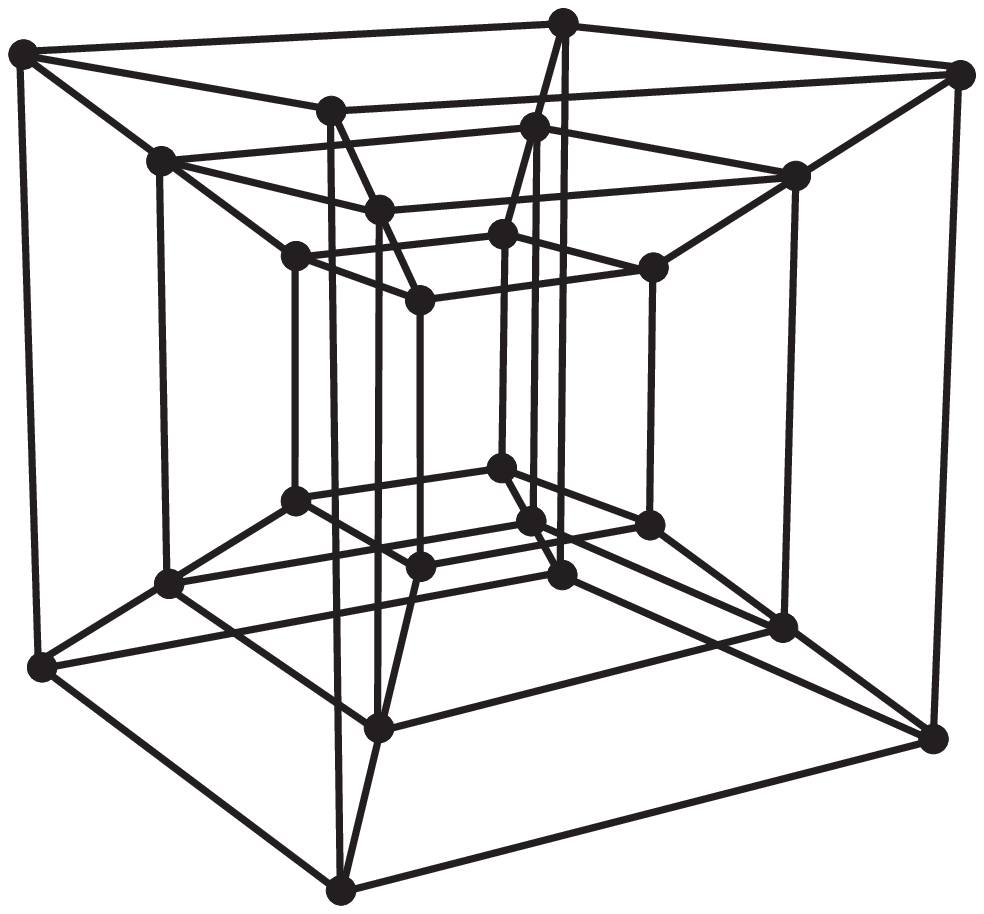}}
\caption{A tesseract-based frame for visualization of the structure of $S_{(4)}$; in order to avoid too crowded appearance of the configuration, there are only shown 24 (out of 81) points and 44 (out  of 108) lines.}
\end{figure}

\pagebreak
In order to find out how many different types of geometric hyperplanes $S_{(4)}$ is endowed with and what the cardinality is of a given hyperplane type, we shall make use of a $3 \times 3 \times 3 \times 3$ grid (`hypercube') representation of $S_{(4)}$, whose rudiments are sketched in Figure 8, and simply blow-up the representatives of both ordinary and extraordinary Veldkamp lines of $S_{(3)}$ whose diagrammatic portrayals were given in the previous section. In particular, taking a representative of the Veldkmap line of $S_{(3)}$ of given type from Figure 6, we first blow-up --- employing the procedure illustrated in Figures 3 and 4 --- all its nine $S_{(2)}$-Veldkamp-sublines
into $S_{(3)}$-hyperplanes to get the character of intersection of the corresponding geometric hyperplane of $S_{(4)}$ with 9 out of its 12 $S_{(3)}$'s, and then look up the character of the remaining three guys in Table 4. Let us illustrate this 
in more detail by taking an explicit example of the type-seven Veldkamp line of $S_{(3)}$, whose representative sits in row four, column one of Figure 6. Let us name its nine $S_{(2)}$-Veldkamp-lines according to the location  of the corresponding $S_{(2)}$'s as left, right, middle-left-right, top, bottom, middle-top-bottom,
front, back, and middle-front-back.  With the help of Figures 3 and 4 we find out that the geometric hyperplanes of $S_{(3)}$ we get by the corresponding blow-ups are of type two, one, four, one, three, three, one, two and four, respectively. From Table 4 we extract that type-seven Veldkamp line of $S_{(3)}$ consists of hyperplanes of type one, three and three. Summarizing, we get that the geometric hyperplane of $S_{(4)}$ we get by blowing-up a representative of Veldkamp line of $S_{(3)}$ of type seven will feature four $S_{(3)}$'s of type $H_1$ and $H_3$, two $S_{(3)}$'s of type $H_2$ and $H_4$, and no $S_{(3)}$'s of type $H_5$ or fully located in the hyperplane. To find the number points and lines this particular geometric hyperplane of $S_{(4)}$ possesses is also quite easy. One simply keeps in mind that each point colored green, blue and red represents a single point of the hyperplane, whilst a point colored black (i.\,e., point of the core) is blown up into three distinct points of the 
hyperplane. As the corresponding subfigure in Figure 6 has nine black points, we find that the hyperplane in question will have $9 \times 3 + (27 - 9) \times 1 = 45$ points. Similarly, if we take into account that each green, blue and red line, as well as each black point, is blown up into a single line of the hyperplane and each black line gives rise to three distinct lines of the hyperplane, we find that the hyperplane must feature 36 lines. We leave it with the interested reader to ascertain the order of individual points of this particular hyperplane type. Following the same procedure, we shall arrive at the same type of geometric hyperplane if we start with a representative of type-eight Veldkamp line of $S_{(3)}$ (shown in row four, column two of Figure 6).
Concerning the question of how many distinct copies of this type of hyperplane we have, it suffices to recall and generalize our observation from Sec.\,3.1 that an ordinary Veldkamp line of $S_{(N-1)}$ gives rise to six different hyperplanes of $S_{(N)}$.
Since our hyperplane, as we have just seen, originates from two different types of ordinary Veldkamp lines of $S_{(3)}$, namely seven and eight, and each of the latter has 324 members (see Table 4), we find that this hyperplane type has as many as $6 \times 2 \times 324 = 3888$ members.
Performing the same analysis with each representative depicted in both Figure 6 and Figure 7, we shall find that geometric hyperplanes of $S_{(4)}$, totaling to 65\,535 (= $2^{16} - 1$ = $|PG(15,2)|$) distinct members, fall into 29 different types whose basic combinatorial properties are collected in Table 5. In the table, the basic ordering is the same as in Tables 1 and 3, i.\,e., the types featuring more points/lines precede those having less ones. However, unlike the previous two cases, we now encounter two or more types having the same number of point/lines, in which case the first goes that having more deep points;
even if this number is the same (for example, types seven and eight), then we list first the one having more points of order three. 
One notes in passing that the hyperplane subject to a detailed examination above belongs to type nine.  
 
Let us make a brief inspection of Table 5 and highlight the most prominent properties of the hyperplane types. Obviously, a type-one hyperplane is a singular hyperplane, that is a hyperplane comprising all the points that are not at maximum distance from a given point; the latter being called the deepest point. On the other hand, a type-29 hyperplane is an ovoid, that is, a set of 27 mutually non-collinear points. Next, let us call, adopting the terminology used for near polygons, a hyperplane of $S_{(4)}$ {\it homogeneous} if all its $S_{(3)}$'s are of the same character. From Table 5 we discern that homogeneous hyperplanes are of types 8, 20 and 29; of course, for the extraordinary hyperplane all of its $S_{(3)}$'s are deep, but there is no homogeneous hyperplane whose $S_{(3)}$'s would be all of type $H_1$ or $H_4$. At the opposite end of the spectrum, there is a type-17 hyperplane, which entails all types of $S_{(3)}$'s except for that corresponding to $S_{(3)}$'s fully located in it. This is also a 
hyperplane type that contains points of all orders; there are other five types enjoying this property, namely 9, 11, 15, 22 and 24. A type-22 hyperplane is also distinguished by the fact that this is the only type whose each member stems from {\it three} different types of Veldkamp lines of $S_{(3)}$. There are 15 types such that each originates from two different types of Veldkamp lines of $S_{(3)}$; each of the remaining 13 types being then a blow-up of a single type of Veldkamp lines of $S_{(3)}$. We further observe that there are only three hyperplane types whose representatives exhibit just a single point of zeroth (i.\,e. minimum) order, and five types with a single point of the fourth (i.\,e. maximum) order. It is also worth mentioning that a hyperplane of type two, three and four, being of weight two, can be gotten as the Veldkamp sum of  any two hyperplanes of type one whose deepest points are at distance two, three and four, respectively.

\begin{table}[t]
\begin{center}
\caption{The types of geometric hyperplanes of $S_{(4)}$ lying on the unique hyperbolic quadric $\mathcal{Q}_0^{+}(15,2) \subset PG(15,2)$ that contains the $S_{(4)}$ (the first orbit) and is invariant under its stabilizer group.} 
{\begin{tabular}{|r|c|c|c|c|c|c|c|c|c|c|c|c|c|r|r|r|c|} \hline \hline
\multicolumn{1}{|c|}{} & \multicolumn{1}{|c|}{} & \multicolumn{1}{|c|}{}  &  \multicolumn{5}{|c|}{}                        & \multicolumn{6}{|c|}{}                      & \multicolumn{1}{|c|}{}  & \multicolumn{1}{|c|}{} & \multicolumn{1}{|c|}{} & \multicolumn{1}{|c|}{} \\
\multicolumn{1}{|c|}{} & \multicolumn{1}{|c|}{} & \multicolumn{1}{|c|}{}  &  \multicolumn{5}{|c|}{$\#$ of Points of Order} & \multicolumn{6}{|c|}{$\#$ of $S_{(3)}$'s of Type} & \multicolumn{1}{|c|}{}  & \multicolumn{1}{|c|}{} & \multicolumn{1}{|c|}{} & \multicolumn{1}{|c|}{} \\
 \cline{4-14}
Tp & Ps & Ls  & 0 & 1 & 2 & 3 & 4 &  D  & $H_1$ & $H_2$ & $H_3$ & $H_4$ &  $H_5$ & VL         & Crd  & BS & W\\
\hline
1  &  65 & 76   &  0  &  0  &  0  &  32 & 33  &  4   &    8   &   0   &   0   &   0   &   0    &   I       & 81   & 2  & 1 \\ 
\hline
2  &  57 & 60   &  0  &  0  & 12  &  24 & 21  &  2   &    6   &   4   &   0   &   0   &   0    &   II,\,1        & 324  & 3  & 2 \\
\hline
3  &  53 & 52   &  0  &  2  & 12  &  26 & 13  &  1   &    6   &   3   &   2   &   0   &   0    &   III,\,2        & 1296 & 4  & 2 \\
\hline 
5  &  49 & 44   &  0  &  8  & 12  &  16 & 13  &  1   &    3   &   6   &   0   &   2   &   0    &   IV,\,4        & 648  & 6  & 3 \\
\hline
7  &  45 & 36   &  0  & 18  &  0  &  18 &  9  &  1   &    0   &   9   &   0   &   0   &   2    &   V,\,6        & 144  & 7  & 3 \\
8  &  45 & 36   &  0  &  0  &  36 &  0  &  9  &  0   &    0   &   12  &   0   &   0   &   0    &  10        & 108  & 17 & 4 \\
9  &  45 & 36   &  2  &  4  &  18 & 16  &  5  &  0   &    4   &   2   &   4   &   2   &   0    &   7,\,8    & 3888 & 8  & 3 \\
10 &  45 & 36   &  0  &  6  &  18 & 18  &  3  &  0   &    3   &   3   &   6   &   0   &   0    &   9,\,11   & 2592 & 9  & 3 \\
\hline
12 &  41 & 28   &  8  &  0  &  24 &  0  &  9  &  0   &    4   &   0   &   0   &   8   &   0    &   14       & 162  & 14 & 4 \\ 
13 &  41 & 28   &  0  &  12 &  18 &  8  &  3  &  0   &    0   &   6   &   4   &   2   &   0    &   15,\,18  & 1944 & 19   & 4 \\ 
14 &  41 & 28   &  0  &  14 &  12 & 14  &  1  &  0   &    1   &   3   &   6   &   2   &   0    &   17,\,21  & 2592 & 15 & 4 \\ 
15 &  41 & 28   &  2  &  8  &  18 & 12  &  1  &  0   &    1   &   3   &   7   &   0   &   1    &   16,\,20  & 2592 & 10   & 3 \\ 
16 &  41 & 28   &  0  &  8  &  24 &  8  &  1  &  0   &    0   &   4   &   8   &   0   &   0    &   19       & 1944 & 20   & 4 \\ 
\hline
21 &  37 & 20   &  4  & 16  & 12  &  0  &  5  &  0   &    0   &   4   &   0   &   8   &   0    &   29       & 972  & 18 & 4 \\ 
22 &  37 & 20   &  4  & 14  & 12  &  6  &  1  &  0   &    0   &   2   &   5   &   4   &   1    &   30--32   & 7776 & 21 & 4 \\ 
23 &  37 & 20   &  3  & 12  & 18  &  4  &  0  &  0   &    0   &   0   &   8   &   4   &   0    &   33       & 3888 & 24 & 4 \\ 
\hline
26 &  33 & 12   &  12 &  12 & 6   &  0  &  3  &  0   &    0   &   2   &   0   &   6   &   4    &   37,\,38  & 648  & 26 & 5 \\
27 &  33 & 12   &  11 &  12 & 6   &  4  &  0  &  0   &    0   &   0   &   4   &   4   &   4    &   39       & 1296 & 29 & 5 \\
\hline \hline
\end{tabular}}
\end{center}
\end{table}

\begin{table}[pth!]
\begin{center}
\caption{Six types of hyperplanes lying on $\mathcal{Q}_{0}^{+}(15,2)$ that in their totality correspond to the image of the set of 2295 maximal subspaces of the symplectic polar space $\mathcal{W}(7,2)$. Interestingly, one orbit consists of homogeneous hyperplanes, viz. of those whose all $S_{(3)}$'s are of type $H_2$.} 
\bigskip
{\begin{tabular}{|r|c|c|c|c|c|c|c|c|c|c|c|c|c|r|r|r|c|} \hline \hline
\multicolumn{1}{|c|}{} & \multicolumn{1}{|c|}{} & \multicolumn{1}{|c|}{}  &  \multicolumn{5}{|c|}{}                        & \multicolumn{6}{|c|}{}                      & \multicolumn{1}{|c|}{}  & \multicolumn{1}{|c|}{} & \multicolumn{1}{|c|}{} & \multicolumn{1}{|c|}{} \\
\multicolumn{1}{|c|}{} & \multicolumn{1}{|c|}{} & \multicolumn{1}{|c|}{}  &  \multicolumn{5}{|c|}{$\#$ of Points of Order} & \multicolumn{6}{|c|}{$\#$ of $S_{(3)}$'s of Type} & \multicolumn{1}{|c|}{}  & \multicolumn{1}{|c|}{} & \multicolumn{1}{|c|}{} & \multicolumn{1}{|c|}{} \\
 \cline{4-14}
Tp & Ps & Ls  & 0 & 1 & 2 & 3 & 4 &  D  & $H_1$ & $H_2$ & $H_3$ & $H_4$ &  $H_5$ & VL         & Crd  & BS & W\\
\hline
1  &  65 & 76   &  0  &  0  &  0  &  32 & 33  &  4   &    8   &   0   &   0   &   0   &   0    &   I        & 81   & 2  & 1 \\ 
\hline
2  &  57 & 60   &  0  &  0  & 12  &  24 & 21  &  2   &    6   &   4   &   0   &   0   &   0    &   II,\,1        & 324  & 3  & 2 \\
\hline
5  &  49 & 44   &  0  &  8  & 12  &  16 & 13  &  1   &    3   &   6   &   0   &   2   &   0    &   IV,\,4        & 648  & 6  & 3 \\
\hline
8  &  45 & 36   &  0  &  0  &  36 &  0  &  9  &  0   &    0   &   12  &   0   &   0   &   0    &  10        & 108  & 17 & 4 \\
\hline
12 &  41 & 28   &  8  &  0  &  24 &  0  &  9  &  0   &    4   &   0   &   0   &   8   &   0    &   14       & 162  & 14 & 4 \\ 
\hline
21 &  37 & 20   &  4  & 16  & 12  &  0  &  5  &  0   &    0   &   4   &   0   &   8   &   0    &   29       & 972  & 18 & 4 \\ 
\hline \hline
\end{tabular}}
\end{center}
\end{table}

A further inspiring insight into the nature of hyperplane types is acquired if one looks at distinguished subconfigurations of the Veldkamp space of $S_{(4)}$, $\mathcal{V}(S_{(4)})$. We have already explicitly seen that  ord-$\mathcal{V}(S_{(N)}) \cong PG(2^N - 1,2)$ for $1 \leq N \leq 3$. It is assumed that this property holds for any $N$, hence ord-$\mathcal{V}(S_{(4)}) \cong PG(15,2)$. In the latter case there, indeed, exists a specific hyperbolic quadric $\mathcal{Q}_0^{+}(15,2) \subset PG(15,2)$ that is composed of 18 particular types of geometric hyperplanes shown in Table 6. Comparing the `VL'- column of this table with the information displayed in Table 4 we find out that these are precisely hyperplane types that originate from those types of (ordinary) Veldkamp lines of $S_{(3)}$ whose cores feature {\it odd} number of points.
This is, however, not the full story as we can look at those orbits of this particular $\mathcal{Q}_0^{+}(15,2)$ that in their totality correspond to the image, furnished by the Lagrangian Grassmannian of type $LGr(4,8)$, of the set of 2295 maximal subspaces of the symplectic polar space $\mathcal{W}(7,2)$ \cite{hsl}. The corresponding hyperplane types are given in Table 7; note that they are exactly the types whose members feature no $S_{(3)}$ of type $H_3$ and $H_5$, each of them also exhibiting points of maximum order (compare with the $S_{(3)}$-case).

\section{Conclusion}
Based on the concept of the (generalized) Veldkamp space of a point-line incidence structure, we have had a detailed look at the structure of the Segre varieties that are $N$-fold direct product of projective lines of size three, $S_{(N)} \equiv PG(1,\,2) \times PG(1,\,2) \times \cdots \times PG(1,\,2)$, for the cases $2 \leq N \leq 4$. In particular, given the fact that $S_{(N)} = PG(1,2) \times S_{(N-1)}$, we have introduced a diagrammatical recipe that shows how to fully recover the properties of Veldkamp points (i.\,e. geometric hyperplanes) of $S_{(N)}$ once we know the types (and cardinalities thereof) of Veldkamp lines of $S_{(N-1)}$. After illustrating this recipe on the $S_{(1)}$-to-$S_{(2)}$ (Sec.\,3.1) and $S_{(2)}$-to-$S_{(3)}$ (Sec.\,3.2) cases, we have made its use to arrive at a complete classification of types of geometric hyperplanes of $S_{(4)}$ (Sec.\,3.3, Table 5). In this latter case, employing the fact that the ordinary part of Veldkamp space of $S_{(N)}$, ord-$\mathcal{V}(S_{(N)})$, is isomorphic to $PG(2^N - 1,2)$, which is conjectured to hold for any $N \geq 1$, 
we  were also able to ascertain (Table 6) which types of geometric hyperplanes of $S_{(4)}$ lie on a particular hyperbolic quadric $\mathcal{Q}_0^+(15,2) \subset PG(15,2) \cong \text{ord-}\mathcal{V}(S_{(4)})$ which contains the $S_{(4)}$  and is invariant under its stabilizer group, as well as to single out those of them (Table 7) that are mapped, via the Lagrangian Grassmannian of type $LG(4,8)$, to the set of 2295 maximal subspaces of the symplectic polar space $\mathcal{W}(7,2)$.

Finally, we would like to stress the following intriguing fact. In our conception of the Veldkamp space of a point-line incidence structure, the extraordinary geometric hyperplane was taken to stand on the same par as any other hyperplane. Otherwise, as it is readily discernible from Tables 1, 3 and 5, by employing our blow-up recipe  it would not only be impossible to recover singular geometric hyperplanes of $S_{(N)}$ from the Veldkamp lines of $S_{(N-1)}$, but, as it is particularly evident in the $N=4$ case (and envisaged to be gradually more pronounced as $N$ increases),  we could not get the correct cardinalities for certain other hyperplane types. The price to be paid for this is that our Veldkamp space is no longer a `homogeneous' geometry, but it comprises a `well-behaving' ordinary part and a rather `ugly' extraordinary part.  This is reminding us distantly of a situation with the concept of an $N$-dimensional projective space over an associative ring with unity, $R$, if one takes the points of this space to be represented by {\it all} free cyclic submodules (FCS's) of the left module $R^{N+1}$; here, the `ordinary' part of the space is generated by unimodular FCS's, whilst its `extraordinary' part --- which is, however, non-empty only for very specific rings, like rings of ternions --- is represented by non-unimodular FCS's (see, e.\,g., \cite{twin,tern} for  some illustrative examples and \cite{cosm} for a possible physics behind).

\section*{Acknowledgment}
This work was partially supported by the VEGA Grant Agency, Project 2/0003/13, as well as by the Austrian Science Fund (Fonds zur F\"orderung der Wissenschaftlichen Forschung (FWF)), Research Project M1564--N27.

\end{document}